\documentclass[12pt, a4paper, reqno]{amsart}
\usepackage{amsmath,amssymb,amsthm,enumitem,mathrsfs,relsize,fontawesome5,fancyhdr}
\usepackage{xcolor,textcomp}
\usepackage{graphicx, epstopdf}
\usepackage{braket,amsfonts}
\usepackage[T1]{fontenc}
\usepackage[utf8]{inputenc}
\usepackage{mathtools}
\usepackage{mwe}
\usepackage{cancel}
\usepackage[top=2.9cm,bottom=2.9cm,left=3cm,right=3cm]{geometry}
\usepackage{hyperref}
\hypersetup{colorlinks=true,linkcolor=red,citecolor=red,filecolor=magenta,urlcolor=blue}
\usepackage[capitalise]{cleveref}
\usepackage{comment}
\usepackage{tikz}
\usepackage{pgfplots}
\usepackage{stmaryrd}
\usepackage[sort,nocompress]{cite}
\numberwithin{equation}{section}

\newcommand{\spann}{\operatorname{span}}

\newcommand{\intt}{\operatorname{int}}
\newcommand{\ri}{\operatorname{ri}}
\newcommand{\co}{\operatorname{co}}
\newcommand{\aff}{\operatorname{aff}}
\newcommand{\rbd}{\operatorname{rbd}}

\newcommand{\Bvee}{\mathlarger{\mathlarger{\vee}}}

\newtheorem{theorem}{Theorem}[section]
\newtheorem{observation}{Observation}[section]
\newtheorem{lemma}[theorem]{Lemma}

\newtheorem{corollary}[theorem]{Corollary}

\newtheorem{example}[theorem]{Example}
 
\title{Differential Inclusions for Gradient and Symmetrized Gradient Operators}
\author[ N. Nesha ]{ Nurun Nesha$^{\dagger,1}$ }
\thanks{$^\dagger$Indian Statistical Institute, Kolkata; $^1$  \faEnvelope \ \href{mailto:nurunnesha55@gmail.com}{nurunnesha55@gmail.com}} 

\pgfplotsset{compat=1.18} 
\begin{document}
\begin{abstract}
 In this article, 
we  study the necessary and sufficient conditions for the  existence of solutions in $W_0^{1,\infty}(\Omega;\mathbb R^n)$ in the minimal dimension of $\spann E$ for the following problem: 
\begin{equation*}
P(D)u\in E \textrm{ a.e. in }\Omega,
\end{equation*}
where $P(D)= D$ or $D+D^{\top}$, and $E\subseteq \mathbb R^{n\times n}$ is a given set. 
We conclude this paper with some properties of  real symmetric matrices.

\end{abstract}
\keywords{differential inclusions, tensor product,  symmetric matrices, symmetric product}
	
\subjclass[2010]{15A69, 15A75, 34A60, 34K09, 34A06}
	
\maketitle
	
\section{Introduction}
\noindent In this article,  we search for a  solution $u\in W_0^{1,\infty}(\Omega;\mathbb R^n)$  of the following differential inclusion:
\begin{equation} \label{19/11/2024 eq1}
    P(D)u\in E \textrm{ a.e. in } \Omega,
\end{equation}
where $\Omega\subseteq \mathbb R^n$ is an open, bounded set, $P(D)$ is the gradient differential operator $D$ or the symmetrized gradient operator $D+D^{\top},$ and $E\subseteq \mathbb R^{n\times n}$ is a given set.

Firstly, we find the minimal dimension of $\spann E$ for the existence of solution of problem \eqref{19/11/2024 eq1}.  Our objective is to determine necessary and sufficient conditions for the existence of solutions of  $\eqref{19/11/2024 eq1}$ on that minimal dimension of $\spann E.$ In \cite{ball1987fine}, Ball--James addresses the inclusion problems for two and three gradient values. Also in \cite{cellina1993minima,cellina1993minimanecessary}, Cellina took into account gradient problems for scalar valued functions and from where we got our motivation to work on the problem \eqref{19/11/2024 eq1}, particularly for the gradient operator. Some results on differential inclusions for gradient operator can also  be seen in Bressan--Flores \cite{BressanFlores}, Dacorogna--Pisante \cite{Pisante_Dacorogna}, Dacorogna--Marcellini \cite{ImplicitPDE}, or Croce\cite{Croce}. 

Though many problems have been studied before for symmetrized gradient operator, the differential inclusions for  symmetrized gradient operator are entirely new and have never been discussed  previously. After adding one constraint that the integral $\int_{\Omega}u\neq 0$, we will see minimal dimension of $\spann E$ and existence of solutions for the problem \eqref{19/11/2024 eq1} at that minimal dimension for this symmetrized gradient operator. The least dimension of $\spann E$ and the existence of a solution at that minimal dimension are the two main topics of our  present study.

The introduction part is covered in the first section and in Section $2$ we will see some of the notations that are used throughout this paper. We have kept the results on gradient and symmetrized gradient operator into two different sections for each one of them. In Section $3$, we study the inclusion problems for gradient operator. The main result in this section states as follows:
\begin{theorem} 
	Let $m,n\in \mathbb N$. Let $\Omega \subseteq \mathbb R^n$ be an open, bounded set$,$ and  $E\subseteq \mathbb R^{m\times n}$  be such that $0\notin E$ with $\dim \spann E=n$. Then$,$ there exists $u\in W_0^{1,\infty}(\Omega;\mathbb R^m)$ satisfying 
    $$D u\in E, \textrm{ a.e. in }\Omega,$$
	if and only if $0\in \ri(\co E)$ and $\spann E=b\otimes \mathbb R^n$ for some $b\in \mathbb R^m \setminus \{0\}$.
\end{theorem}
In Section $4$, we consider the symmetrized gradient operator $D+D^{\top}$ and the main theorem of this section is as followw: 
\begin{theorem} Let $n\in \mathbb N$, $\Omega\subseteq \mathbb R^n$ an open, bounded set, and $E\subseteq \Bvee^2(\mathbb R^n)\setminus \{0\}$ be such that $\dim \spann E=n$. There exists a solution $u\in W_0^{1,\infty}(\Omega;\mathbb R^n)$ such that 
	\begin{align}
	&Du+(Du)^{\top}\in E,\\
    & \textrm{ and } \int_{\Omega}u\neq 0 \nonumber
	\end{align}
	if and only if $0\in\ri(\co E)$ and $(\spann E)^{\perp}b=0$ for some $b\in \mathbb R^n \setminus \{0\}.$
\end{theorem}

\noindent We conclude this paper with few observations of real symmetric matrices in our final appendix section.

The results of differential inclusion problems can be applied, embracing the notions due to Cellina \cite{cellina1993minima,cellina1993minimanecessary} and Friesecke \cite{Friesecke}, to obtain solutions for a non-convex variational problem. 
\noindent In variational analysis, differential inclusions come up in the following way. Let us suppose that we are interested on the existence of minimizer $u$ of the following problem 
$$
(P)\ \ \ \ \left\{\int_{\Omega}f(Du):u\in W^{1,\infty}_{0}(\Omega;\mathbb{R}^m)\right\},
$$ 
where $f:\mathbb{R}^{m\times n}\to\mathbb{R}$ is locally bounded, lower semicontinuous and non-negative. If $f$ is quasiconvex, then $(P)$ has clearly a minimizer, namely, the zero function. The problem becomes much more interesting when $f$ is not quasiconvex. When $f$ is not quasiconvex, it follows from Relaxation Theorem that $(P)$ has a minimizer if and only if there exists  $\bar{u}\in W^{1,\infty}_{0}(\Omega;\mathbb{R}^m)$ satisfying 
$$
f(D\bar{u})=Qf(D\bar{u}),\text{ a.e. in }\Omega,
$$
and $$
\int_{\Omega}Qf(D\bar{u})=Qf(0). \mathcal L^n(\Omega).
$$
In other words, to ensure the existence of a minimizer, we have to study the differential inclusion that we want  $\bar{u}\in W^{1,\infty}_{0}(\Omega;\mathbb{R}^m)$ satisfying 
$$D\bar{u}\in K, \textrm{ where }K:=\{\xi\in \mathbb{R}^{m\times n}: Qf(\xi)=f(\xi)\}.$$
More details can be seen in Dacorogna \cite{dacorogna2005nonconvex}.
\noindent Many results as well as applications of differential inclusions can be observed in Dacorogna--Pisante \cite{PisanteDacorogna}, Dacorogna--Fonseca \cite{DacorognaFonseca}, Dacorogna--Marcellini \cite{DacorognaMarcellini}, Blasi--Pianigiani \cite{BlasiPianigiani}, Sil \cite{SwarnenduACV}, and Sychev \cite{Sychev}.

\section{Notations}
\noindent We gather here some notations which will be used throughout this article. For more details on exterior algebra and differential forms see \cite{csato2011pullback} and for convex analysis see \cite{DirectBernard} or \cite{Rockafellar}.
\begin{itemize}
\item [$\bullet$] For subsets $C,F\subseteq \mathbb R^n$, $\co C $, $\aff C$ denote the convex hull and affine hull of $C$ respectively, and $\intt_{F}C$ denotes  the interior of $C$ with respect to the topology relative to $F.$  
\item [$\bullet$] For a  convex set  $C\subseteq \mathbb R^n$, $\ri( C)$ denotes the relative interior of $C$ which is the interior of $C$ with respect to the topology relative to affine hull of $C.$ Equivalently, $\ri C=\intt_{\aff C}C$. Similarly, $\rbd ( C)$ means relative boundary of $C$.
    \item [$\bullet$] $a\otimes b$ means $ab^{\top}$ and $ab^{\top}\in \mathbb R^{n\times m}$ for $a\in \mathbb R^n,$ $b\in \mathbb R^m.$
    \item [$\bullet$] $a\vee b:= a\otimes b+b \otimes a$ for $a,b \in \mathbb R^n$ and reads as $a$ symmetric product $b.$
    \item [$\bullet$] $a\wedge b:=a\otimes b-b \otimes a$ for $a,b\in \mathbb R^n$ and reads as $a $ wedge product $b.$
    \item [$\bullet$]  $\lrcorner$ and $\langle;\rangle$  denote the  interior product and the scalar product  respectively.
    \item [$\bullet$] $\mathbb R^n \vee b=\{x\vee b: x\in \mathbb R^n\}$, $\mathbb R^n \wedge b=\{x\wedge b: x\in \mathbb R^n\}$.
    \item [$\bullet$] $\mathbb R^{n\times n}_{\textrm{sym}}$ denotes as space of all $n\times n$ real symmetric matrices.
    \item [$\bullet$]  $GL_n^{\textrm{sym}}(\mathbb R)=\{A\in \mathbb R_{\textrm{sym}}^{n\times n}: \det A\neq 0\}.$
    \item [$\bullet$] $\Bvee^2(\mathbb R^n)=\spann \{e^i\vee e^j:1\leq i\leq j \leq n\}$, where $\{e^i:i=1,\ldots , n\}$ is the canonical basis of $\mathbb R^n.$ So  $\Bvee^2(\mathbb R^n)\cong\mathbb R^{n\times n}_{\textrm{sym}}.$ 
    \item [$\bullet$] For $u:\mathbb R^n \to \mathbb R^m$, $Du$ denotes the gradient  of $u$.
    \item $\mathcal L^n(A)$ is the Lebesgue measure of $A\subseteq \mathbb R^n.$
    \item [$\bullet$] $u\in W^{1,\infty}_{0}(\Omega)$ means $u: \Omega \to \mathbb R$ and $W_0^{1,\infty}$ denotes the standard Sobolev space as defined in \cite{EvansPDE}.
    \item  [$\bullet$] For $A\subsetneq B$, it means that $A$ is a proper subset of $B.$
\end{itemize}

\section{Gradient Operator}
\noindent In this section, we obtain  necessary and sufficient conditions for the existence of solutions $u\in W_0^{1,\infty}(\Omega;\mathbb R^m)$ to the following differential inclusion problem: 
\begin{equation*}
Du \in E,
\textrm{ a.e. in }\Omega,
\end{equation*}
where $E \subseteq \mathbb R^{m\times n}\setminus \{0\}$ is a given set and $\Omega \subseteq \mathbb R^n$ is an open, bounded set. We begin with the following elementary result as a lemma.
\begin{lemma} \label{8/3/2024 lemma 1}
    Let $m,n,p\in \mathbb N$. Let $\{a_1,\ldots, a_p\}\subseteq \mathbb R^n$ be a linearly independent set and $\{b_1,\ldots, b_p\}\subseteq \mathbb R^m \setminus \{0\}.$ Then $\{b_i\otimes a_i :i=1,\ldots,p \}\subseteq \mathbb R^{m\times n}$ is linearly independent.
\end{lemma}
\begin{proof}
    Let us suppose that 
   \begin{equation}\label{Nov 12 2024 eq1}
   \displaystyle \sum_{i=1}^{p}\lambda_i b_i \otimes a_i=0 \textrm{ and } b_i=\displaystyle \sum_{j=1}^{m}b_{i}^{j}e^{j},
    \end{equation}
    where $\lambda_i\in \mathbb R $ for $i=1,\ldots , p$ and $b_{i}^{j}\in \mathbb R$ for $i=1,\ldots,p$; $j=1,\ldots,m$. If we rearrange the above equation \eqref{Nov 12 2024 eq1}, we can state that  
    \begin{equation}\label{Nov 12 2024 eq2}
    \displaystyle\sum_{i=1}^{p}\lambda_i (a_i\otimes b_i) e^j=0 \textrm{ i.e., }\displaystyle\sum_{i=1}^{p}\lambda_i a_i \langle b_i;e^j \rangle =0,
    \end{equation}
    for all $j=1,\ldots,m$, 
    where $\{e^j: j=1,\ldots, m\}$ is the canonical basis of $\mathbb R^m$. Consequently, the above equations $\eqref{Nov 12 2024 eq1}$ and \eqref{Nov 12 2024 eq2} provide us with   $\lambda_i b_i^{j}=0$ for all $i=1,\ldots, p$; $j=1,\ldots,m$ as $\{a_1,\ldots, a_p\}$ is a linearly independent set. Since $\{b_1,\ldots, b_p\}\subseteq \mathbb R^m \setminus \{0\}$, we can conclude that  $\lambda_i=0$ for all $i=1,\ldots,p.$ Hence $\{b_i\otimes a_i: i=1,\ldots,p\}$ is linearly independent.
\end{proof}
\begin{corollary} \label{13/11/2024 cor 1}
	Let $m,n,p\in \mathbb N$. Let $\{a_1,\ldots, a_p\}\subseteq \mathbb R^n$ be a linearly independent set and $\{b_1,\ldots, b_p\}\subseteq \mathbb R^m$. If $\ \displaystyle \sum_{i=1}^{p}b_i\otimes a_i=0,$ then  $b_i=0$ for all $i=1,\ldots, p.$ 
\end{corollary}
\noindent Next we will present two fundamental  lemmas \cite[Lemma 3.1 and 3.2]{bandyopadhyay2015some} with their obvious proofs and then we will prove one basic lemma in  \Cref{20/05/2024 lemma 1} below which was originated by a result  proved in \cite{bandyopadhyay2015some} (cf. Lemma $3.3$) for differential forms.
\begin{lemma} \label{21/05/2024 lemma 1}
	Let $X$ be a normed linear space, $Y\subseteq X$ be a linear subspace, and let $A\subseteq X.$ Then$,$ $0\in \intt_{Y}\co (A\cap Y)$ if and only if $$\spann (A\cap Y)=Y, \textrm{ and } 0\in \ri\co (A\cap Y).$$
	
\end{lemma}
\begin{lemma}[\textbf{Carath\'eodory}] \label{Carath\'eodory 22/06/24}
	Let $X$ be a normed linear space, and $E\subseteq X.$ Then$,$ $0\in \ri (\co E)$ if and only if there exist $m\in \mathbb N,$ $m\geq \dim \spann E+1,$ $z^i\in E,$ $t^i >0$ for every $1\leq i\leq m,$ such that $$\displaystyle \sum_{i=1}^{m}t^iz^i=0, \ \sum_{i=1}^{m}t^i=1, \textrm{ and } \spann \{z^1,\ldots, z^m\}=\spann E.$$
\end{lemma}
\noindent We will apply these two lemmas mentioned before to present the lemma below.
\begin{lemma} \label{20/05/2024 lemma 1}
	Let $m,n \in \mathbb N,$ $b\in \mathbb R^m\setminus \{0\},$ and $E\subseteq \mathbb R^{m\times n}$.  Assume  the following condition holds:  $$0\in \intt_{b\otimes \mathbb R^n} \co [E\cap (b\otimes \mathbb R^n)].$$  
	Then there exists $F\subseteq \mathbb R^n$ such that $$E\cap (b\otimes \mathbb R^n)=b\otimes F, \textrm{ and } 0\in \intt \co F.$$
\end{lemma}
\begin{proof}
Let us take $F:=\{x\in \mathbb R^n:b\otimes x \in E\}$. Clearly, $E\cap (b\otimes \mathbb R^n)=b\otimes F$. As $0\in \intt_{b\otimes \mathbb R^n}\co \left[ E\cap (b\otimes \mathbb R^n) \right],$ making use of  \Cref{21/05/2024 lemma 1} above, we can assert  that 
	\begin{equation*}
	\spann [E\cap (b\otimes \mathbb R^n)]=b\otimes \mathbb R^n,
	\textrm{ and } 0\in \ri \co [E\cap (b\otimes \mathbb R^n)]=\ri \co [b\otimes F].
	\end{equation*}
	In turn$,$ employing  \Cref{Carath\'eodory 22/06/24}, we identify $m\in \mathbb N,$ $m\geq \dim \spann [b\otimes F]+1,$ $\lambda_i>0,$ and $f_i \in F$ for  $i=1,\ldots, m$ so that 
	\begin{align}\label{Nov 13 2024 eq 1}
	& \displaystyle \sum_{i=1}^{m}\lambda_i =1, \ \ \sum_{i=1}^{m}\lambda_i (b\otimes f_i)=0,\\
	\textrm{ and } & \spann \{b\otimes f_i:i=1,\ldots, m\}=\spann \left [b\otimes F\right ]. \nonumber
	\end{align}
	To establish that $0\in \ri\co F,$ thanks to the \Cref{Carath\'eodory 22/06/24} as it suffices to show the following: 
    $$m\geq \dim \spann F +1, \ \displaystyle \sum_{i=1}^{m}\lambda_i f_i =0,\textrm{ and }\spann \{f_i: i=1,\ldots, m\}=\spann F.$$
    As the linear map $T: \mathbb R^n \to \mathbb R^{m\times n}$ defined as $T(x):=b\otimes x,$ for all $x\in \mathbb R^n$ is one-one, this means $$\dim \spann F =\dim \spann [b\otimes F].$$
	Hence$,$ $m\geq \dim \spann F+1$. Additionally, we have $\displaystyle \sum_{i=1}^{m}\lambda_i f_i=0$ through the equation \eqref{Nov 13 2024 eq 1}. It remains to check that $\spann \{f_i: i=1,\ldots, m\}=\spann F$.
	The fact that the following is reasonable:
	\begin{align*}
	b\otimes \spann \{f_i: i=1,\ldots,m\}&=\spann \{b\otimes f_i:i=1,\ldots, m\}\\
	& =\spann [b\otimes F]=b\otimes \spann F,
	\end{align*}
	we can write  $\spann \{f_i: i=1,\ldots, m\}=\spann F$. Now we have from above that  $b\otimes \mathbb R^n=\spann [b\otimes F]=b\otimes \spann F.$ Consequently, $\spann F=\mathbb R^n$ and \Cref{21/05/2024 lemma 1} gives us $0\in \intt \co F$. 
    This concludes the proof.
\end{proof}
\noindent Regarding the scalar gradient situation, we recall a result below (see Lemma $2.11$ in \cite{ImplicitPDE} for proof).
\begin{lemma}\label{25/11/2024 lemma1}
    Let $n\in \mathbb N$, $\Omega \subseteq \mathbb R^n$ be an open, bounded set, and $F\subseteq \mathbb R^n$. Then there exists $u\in W_0^{1,\infty}(\Omega)$ such that $$Du \in F\textrm{ a.e. in }\Omega$$
    if and only if $0\in F\cup \intt \co F.$  Additionally, $u$ can be taken piecewise affine, $u\geq 0$, and $\int_{\Omega}u>0$ if $0\in \intt \co F$.
\end{lemma}

\noindent We will now present the first key result of this section. 
\begin{theorem}[\textbf{Minimum dimension}] \label{02/4/2024 theorem 1}
	Let $m,n\in \mathbb N.$ Let   $E\subseteq \mathbb R^{m\times n}\setminus \{0\}$ and $\Omega \subseteq \mathbb R^n$ be an open, bounded set. If there exists $u\in W_0^{1,\infty}(\Omega , \mathbb R^m)$ such that $$Du\in E \textrm{ a.e. in } \Omega,$$ then $\dim \spann E\geq n.$ Moreover, if we take $\dim \spann E=n,$ then there exists $b\in \mathbb R^m \setminus \{0\}$ such that $\spann E=b \otimes \mathbb R^n.$
\end{theorem}
\begin{proof}[Proof of \Cref{02/4/2024 theorem 1}] 
	Let $\mathbb P: \mathbb R^{m\times n}\to \mathbb R^{m\times n}$ denote the projection onto the orthogonal complement of $\spann E.$ Since $u\in W_{0}^{1,\infty}(\Omega; \mathbb R^m)$, extending $u$ by $0$ to $\mathbb R^m$ we can write $$\mathbb P \left(Du(x)\right)=0 \textrm{ a.e. in } \mathbb R^n.$$
	Applying Fourier transform, we get $$\mathbb P( \hat{u}(x)\otimes x)=0 \textrm{ for a.e. } x\in \mathbb R^n,$$
	$$\textrm{i.e., } \hat{u}(x) \otimes x\in \spann E \textrm{ a.e. in }\mathbb R^n.$$
	Let $\hat{u}(x_0)\neq 0$ for some $x_0 \in \mathbb R^n.$ Then there exists a neighbourhood of $x_0$ and hence a basis $\{x_1,\ldots, x_n\}$ of $\mathbb R^n$  such that $\hat u(x_i)\neq 0$ for all $i=1,\ldots,n.$ Therefore from  \Cref{8/3/2024 lemma 1} we can claim that $\{ \hat{u}(x_i)\otimes x_i: i=1,\ldots, n\}$ is linearly independent and hence $\dim \spann E\geq n.$
	
\noindent	Now let $\dim \spann E=n.$ Then we will show that there exists $b\in \mathbb R^m\setminus \{0\}$ such that $\spann E=b\otimes \mathbb R^n.$ In regard to do this, let $\mathcal S\subseteq \mathbb R^n$ be such that $\mathcal L^n(\mathcal S)>0$ and $\hat{u}(x)\otimes x\in \spann E \setminus\{0\}$ for all $x\in \mathcal S.$ As $\mathcal L^n(\mathcal S)>0$, there exists a basis $\{x_1,\ldots,x_n\}$ of $\mathbb R^n$ such that $\hat{u}(x_i)\otimes x_i\in \spann E\setminus\{0\}$ and $\{\hat{u}(x_i)\otimes x_i:i=1,\ldots,n\}$ is a basis of $\spann E$ (utilizing  \Cref{8/3/2024 lemma 1}). For each $J\subsetneq \{1,\ldots,n\}$, we identify 
$$E_J:=\spann \{x_j:j\in J\}.$$
Evidently, $\mathcal L^n{(E_J)}=0$  for all $J\subsetneq \{1,\ldots,n\}$ and also $\mathcal L^n(\mathcal S)>0$. Thus we are able to write
$$\mathcal L^n\bigg(\mathcal S\setminus \displaystyle\bigg(\bigcup\limits_{J\subsetneq \{1,\ldots,n\}}E_J\bigg)\bigg)>0.$$
Let us choose $y\in \mathcal S\setminus\displaystyle\bigcup\limits_{J\subsetneq \{1,\ldots,n\}}E_J\ $ and $y=\displaystyle\sum_{i=1}^{n}\alpha_i x_i$, where $\alpha_i\in \mathbb R$ for all $i=1,\ldots,n.$ Then $\alpha_i \neq 0$ for all $i=1,\ldots,n$ and $\hat{u}(y)\otimes y\in \spann E\setminus \{0\}.$ Let 
$$\hat{u}(y)\otimes y=\displaystyle\sum_{i=1}^{n}\lambda_i\hat{u}(x_i)\otimes x_i,$$
where $\lambda_i\in \mathbb R$ for $i=1,\ldots,n$. This provides  us  $\hat{u}(y)\otimes \Big(\displaystyle\sum_{i=1}^{n}\alpha_i x_i\Big)=\displaystyle\sum_{i=1}^{n}\lambda_i\hat{u}(x_i)\otimes x_i,$ i.e., $\displaystyle\sum_{i=1}^{n}\left(\alpha_i \hat{u}(y)-\lambda_i \hat{u}(x_i)\right)\otimes x_i=0.$ Implementing   \Cref{13/11/2024 cor 1}, we can state  that $\alpha_i \hat{u}(y)-\lambda_i \hat{u}(x_i)=0$ for all $i=1,\ldots,n$.   Consequently, $$\hat{u}(x_i)=\frac{\alpha_i}{\lambda_i}\hat{u}(y) \textrm{ for all }i=1,\ldots,n.$$
	Hence $\spann E=\spann \{\hat{u}(x_i)\otimes x_i:i=1,\ldots,n\}=\spann \{\hat{u}(y)\otimes x_i:i=1,\ldots,n\}=\hat{u}(y)\otimes \mathbb R^n.$ This finishes the proof. 
\end{proof}

\noindent We are now all set to demonstrate the core theorem of this section.
\begin{theorem}[\textbf{Existence theorem in the minimal dimension}] \label{20/05/2024 thm 1'} 
	Let $m,n\in \mathbb N$. Let $\Omega \subseteq \mathbb R^n$ be an open, bounded set$,$ and  $E\subseteq \mathbb R^{m\times n}$  be such that $0\notin E,$ and $\dim \spann E=n$. Then$,$ there exists $u\in W_0^{1,\infty}(\Omega;\mathbb R^m),$ satisfying $$D u\in E, \textrm{ a.e. in }\Omega,$$
	if and only if $0\in \ri(\co E)$ and $\spann E=b\otimes \mathbb R^n$ for some $b\in \mathbb R^m \setminus \{0\}$.
\end{theorem}
\begin{proof}
	When there exists $u\in W_0^{1,\infty}(\Omega;\mathbb R^m),$ satisfying 
    \begin{equation} \label{Nov 14 2024 eq1}
    Du \in E \textrm{ a.e. in }\Omega, 
    \end{equation}
    and $\dim \spann E=n$, then the conclusions that there exists $b\in \mathbb R^m \setminus \{0\}$ and $\spann E=b\otimes \mathbb R^n,$ are derived from  \Cref{02/4/2024 theorem 1}. Now$,$ let $f: \mathbb R^{m\times n}\to \mathbb R$ be a convex function such that $f|_{E}\leq 0$. Afterwards$,$ making use of  Jensen's inequality and the inclusion \eqref{Nov 14 2024 eq1}, we are able to write 
    $$f(0)=f\bigg(\displaystyle \int_{\Omega} Du(x) dx\bigg)\leq \int_{\Omega}f\left(Du\right)(x) dx\leq 0.$$
    As a result, we obtain $f(0)\leq 0$. Together with this  and the Proposition $2.36$ of Dacorogna \cite{DirectBernard} allow  us to write $0\in \overline{\co E}$.
	
\noindent	Now, let us suppose that $0\notin \ri(\co E),$ then $0\in \rbd (\co E)$. Next, taking advantage of Separation Theorem (cf. Theorem $2.10$ in \cite{DirectBernard}), there exists $P\in \spann E\setminus \{0\}$ such that $$\langle A; P \rangle\geq 0, \textrm{ for all }A\in \overline{\co E}.$$
	In particular, $\langle Du(x); P \rangle \geq 0,$ for a.e. $x\in \Omega$. However, as $u\in W^{1,\infty}_{0}(\Omega; \mathbb R^m),$ we are able to  write 
    $$\displaystyle \int_{\Omega}\langle Du(x);P\rangle=0.$$
	Thereby, $\langle Du(x);P\rangle=0, \textrm{ for a.e. }x\in \Omega.$
 \noindent This suggests that $P\in (\spann E)^{\perp},$ which is a contradiction, as $P\in \spann E \setminus \{0\}.$ Consequently$,$ $0\in \ri (\co E)$.
	
   Conversely, let us assume that $0\in \ri (\co E)$ and $\spann E= b\otimes \mathbb R^n$ for some $b\in \mathbb R^m \setminus \{0\}$. Then$,$ applying  \Cref{21/05/2024 lemma 1}, we can write $0\in \intt_{b\otimes \mathbb R^n}\co \left( E\cap (b\otimes \mathbb R^n) \right)$. Following that, we can employ  \Cref{20/05/2024 lemma 1} and we identify   $F\subseteq \mathbb R^n$ such that $E=E\cap (b\otimes \mathbb R^n)=b\otimes F$ as well as  $0\in \intt \co F$.  
 And finally, invoking \Cref{25/11/2024 lemma1}, there exists $v\in W^{1,\infty}_{0}(\Omega),$ satisfying $Dv \in F,$ a.e. in $\Omega$. Let us define $u\in W^{1,\infty}_{0}(\Omega;\mathbb R^m)$ by 
 \begin{equation*}
 u(x):=v(x)b.
 \end{equation*}
Immediately$,$ $D u\in E,$ a.e. in $\Omega$. This concludes the proof.
\end{proof}
\section{Symmetrized Gradient Operator}
\noindent In this section, we will   introduce differential inclusions for  symmetrized gradient operator. Basically, we are interested in finding  $u\in W_0^{1,\infty}(\Omega;\mathbb R^n)$ for the following differential inclusion problem:
\begin{equation}\label{01/12/2024 eq1}
    Du+(Du)^{\top}\in E, \textrm{ a.e. in }\Omega \textrm{ and }
    \int_{\Omega}u\neq 0,
\end{equation}
where $E\subseteq \Bvee^2(\mathbb R^n)\setminus \{0\}$ and $\Omega\subseteq \mathbb R^n$ being an open, bounded set. The reason behind taking the constraint  $\int_{\Omega}u\neq 0$ can be seen in Section $2$ of \cite{bandyopadhyay2015some}. Let us begin with the following elementary lemma without proof. 
\begin{lemma} \label{15/11/2024 lemma 1}
	Let $n\in \mathbb N,$ $b\in \mathbb R^n\setminus \{0\}$. Then $\dim (\mathbb R^n \vee b)=n.$
\end{lemma}
   
\noindent In the following lemma, we will prove necessary and sufficient condition for an $n$ dimensional subspace $S\subseteq \mathbb R^{n\times n}_{\textrm{sym}}$ to be of the form $\mathbb R^n \vee b$. 
\begin{lemma} \label{17/11/2023 lemma 1}
Let $n\geq 2$, $b\in \mathbb R^n \setminus \{0\},$ and $S$ be an $n$ dimensional subspace of $\mathbb R^{n\times n}_{\textrm{sym}}$. Then $S=\mathbb R^n \vee b$ if and only if $S^{\perp}b=0$ i.e., $b\in \ker A$ for all $A\in S^{\perp}.$
\end{lemma}
\begin{proof}
Let us assume that $S^{\perp}b=0$. For $A\in S^{\perp}$ and  $x\in \mathbb R^n$, we can write   $\langle A;x \vee b \rangle=0$, because $\langle A;x \vee b \rangle =\langle Ab;x\rangle +\langle Ax;b \rangle=2\langle Ab;x \rangle , \textrm{ as } A^{T}=A$.
Consequently, $A\in \left(\mathbb R^n\vee b\right)^{\perp}$, i.e., $S^{\perp}\subseteq (\mathbb R^n \vee b)^{\perp}.$ This indicates that $\mathbb R^n \vee b \subseteq S.$ Thus employing  \Cref{15/11/2024 lemma 1},  we obtain $\mathbb R^n \vee b =S.$ 

 Conversely, if $S=\mathbb R^n \vee b$, then for $A\in S^{\perp}$ and  $x\in \mathbb R^n$, we can say likewise $\langle A;x\vee b \rangle =0$, i.e., $\langle Ab;x \rangle =0$ and this is true for all $x\in \mathbb R^n$. Therefore, $Ab=0$. This finishes our proof.
\end{proof}
\noindent We will note down  a lemma below  as similar to  \Cref{20/05/2024 lemma 1} for symmetric product. Since the proof is identical to \Cref{20/05/2024 lemma 1} if we substitute the symmetric product for the tensor product, we write it down without proof.
\begin{lemma}\label{18/11/2024 lemma1}
Let $n \in \mathbb N,$ $b\in \mathbb R^n\setminus \{0\},$ and $E\subseteq \mathbb R^{n\times n}$.  Assume  the following condition holds:  $$0\in \intt_{b\vee \mathbb R^n} \co [E\cap (b\vee \mathbb R^n)].$$  
	Then there exists $F\subseteq \mathbb R^n$ such that 
    $$E\cap (b\vee \mathbb R^n)=b\vee F \textrm{ and } 0\in \intt \co F.$$
\end{lemma}
\noindent In the following theorem, we will prove minimal dimension of $\spann E$ for the existence of solution of the differential inclusion problem \eqref{01/12/2024 eq1} adding one constraint that $\int\limits_{\Omega}u\neq 0.$ 
\begin{theorem}[\textbf{Minimum dimension with constraint}]\label{24/06/2024 thm1}
	Let $n\in \mathbb N,$ $\Omega\subseteq \mathbb R^n$ be an open, bounded set$,$ and $E$ being a subset of $\Bvee^{2}(\mathbb R^n)$ such that $0\notin E.$ Let  $u\in W_0^{1,\infty}(\Omega; \mathbb R^n)$ be a solution of the following problem$,$ 
	\begin{align*}
		& Du+(Du)^{\top}\in E \textrm{ a.e. in }\Omega, \\
		\textrm{and } & \int_{\Omega} u \neq 0.
	\end{align*}
Then $\dim \spann E\geq n$. 
\end{theorem}
\begin{proof}
	Let $\mathbb P:\Bvee^2(\mathbb R^n)\to \Bvee^2(\mathbb R^n)$ be the projection onto the orthogonal complement of $\spann E$. Since $u \in W_{0}^{1,\infty}(\Omega;\Bvee^2)$, extending $u$ by $0$ to $\mathbb R^n,$ it follows that 
    $$\mathbb P\left(Du +(Du)^\top\right)=0, \textrm{ a.e. in } \mathbb R^n.$$
    Now$,$ let us apply the Fourier transformation. After that$,$ we get 
    $$\mathbb P\left( x\vee \left[\int_{\mathbb R^n}u(y)\cos \left(2\pi\langle x;y\rangle\right) \mathrm{d}y\right]\right)=0, \textrm{ for all }x\in \mathbb R^n,$$
    which is equivalent to $$x\vee \left[\int_{\mathbb R^n}u(y)\cos \left(2\pi\langle x;y \rangle\right) \mathrm{d}y\right]\in \spann E ,\textrm{ for every } x\in \mathbb R^n.$$ 
	Let $$f(x)=\int_{\mathbb R^n}u(y)\cos \left(2\pi \langle x;y \rangle\right) \mathrm{d}y.$$ From here, we get $$f(0)=\int_{\Omega}u \neq 0.$$
	Now$,$ we will show that 
	\begin{align}\label{20/11/2023 12:55 pm}
		\mathbb R^n \vee f(0)\subseteq \spann \{x\vee f(x):x\in \mathbb R^n\}.
	\end{align}
	To prove \eqref{20/11/2023 12:55 pm}, let $x\in \mathbb R^n \setminus \{0\}$ be fixed. Note that$,$ for every $\lambda \in \mathbb R\setminus \{0\},$ $$x\vee f(\lambda x)=\frac{1}{\lambda}\left[\lambda x \vee f(\lambda x) \right]\in \spann \{y\vee f(y): y\in \mathbb R^n\}.$$ 
	Since $\spann \{y\vee f(y):y\in \mathbb R^n\}$ is closed and $f$ is continuous at $0$, it follows that $$x\vee f(0)\in \spann \{y\vee f(y):y\in \mathbb R^n\},$$ 
	after letting $\lambda \to 0$.
	
\noindent	Therefore$,$ 
$$\mathbb R^n \vee f(0)\subseteq \spann \{y\vee f(y):y\in \mathbb R^n\}.$$
Hence$,$ our assertion \eqref{20/11/2023 12:55 pm} is proved. Now$,$ as we know  $\dim (\mathbb R^n \vee f(0)=n,$ it follows that $\dim \spann E\geq n$.
\end{proof}
\noindent Next, we are going to show necessary and sufficient conditions for  the existence of solutions of the symmetrized gradient operator in the minimal dimension.
\begin{theorem}[\textbf{Existence of solution for symmetrized gradient}] Let $n\in \mathbb N$, $\Omega\subseteq \mathbb R^n$ an open, bounded set, and $E\subseteq \Bvee^2(\mathbb R^n)\setminus \{0\}$ be such that $\dim \spann E=n$. There exists a solution $u\in W_0^{1,\infty}(\Omega;\mathbb R^n)$ such that 
	\begin{align}\label{17/11/2024 eq1}
	&Du+(Du)^{\top}\in E,\\
    & \textrm{ and } \int_{\Omega}u\neq 0 \nonumber
	\end{align}
	if and only if $0\in\ri(\co E)$ and $(\spann E)^{\perp}b=0$ for some $b\in \mathbb R^n \setminus \{0\}.$
\end{theorem}
\begin{proof}
    Utilizing  \Cref{24/06/2024 thm1}, we get $\spann E=\mathbb R^n \vee b$, where $b=\int_{\Omega}u\in \mathbb R^n \setminus \{0\}$. Consequently,  \Cref{17/11/2023 lemma 1} gives us $(\spann E)^{\perp}b=0$. Now let $f:\mathbb R^{n\times n}\to \mathbb R$ be a convex function such that $f|_{E}\leq 0$. Then applying Jensen's inequality and the inclusion \eqref{17/11/2024 eq1}, we can write 
    \begin{equation*}
        f(0)=f\left(\int_{\Omega}\left(Du+(Du)^{\top}\right)(x)dx\right) \leq \int_{\Omega} f(Du+(Du)^{\top})(x) dx \leq 0.
    \end{equation*}
    Thus, we obtain $f(0)\leq 0$ and hence $0\in \overline{\co E}$ following Proposition $2.36$ of Dacorogna \cite{DirectBernard}. 
    
    \noindent Now, let us suppose that $0\notin \ri(\co E),$ then $0\in \rbd (\co E)$. Next, taking advantage of separation Theorem (cf. Theorem $2.10$ in \cite{DirectBernard}), there exists $P\in \spann E\setminus \{0\}$ such that 
    $$\langle A; P \rangle\geq 0 \textrm{ for all }A\in \overline{\co E}.$$
	In particular, $\langle Du+(Du)^{\top}(x); P \rangle \geq 0$ for a.e. $x\in \Omega$. However, as $u\in W^{1,\infty}_{0}(\Omega; \mathbb R^m),$ we are able to  write 
    $$\displaystyle \int_{\Omega}\langle \left(Du+(Du)^{\top}\right)(x);P\rangle=0.$$
	Thereby, $\langle Du+(Du)^{\top}(x);P\rangle=0 \textrm{ for a.e. }x\in \Omega.$
 \noindent This suggests that $P\in (\spann E)^{\perp},$ which is a contradiction, as $P\in \spann E \setminus \{0\}.$ Consequently$,$ $0\in \ri (\co E)$.

 Conversely, if $0\in \ri (\co E)$ and $(\spann E)^{\perp} b=0$ for some $b\in \mathbb R^n \setminus \{0\}$, then applying  \Cref{18/11/2024 lemma1} and  \Cref{17/11/2023 lemma 1}, we can say that there exists $F\subseteq \mathbb R^n$ such that $0\in \intt \co F$ and $E\cap(b \vee \mathbb R^n )=b \vee F$. Now, calling forth \Cref{25/11/2024 lemma1}, we get a $v\in W_0^{1,\infty}(\Omega)$ such that $Dv \in F$ a.e. in $\Omega$. Let us define $u:=vb\in W^{1,\infty}_0(\Omega;\mathbb R^n)$. Hence $Du+(Du)^{\top}\in E$ a.e. in $\Omega$. This concludes the proof.
\end{proof}
\section{Appendix}
\noindent In this section, we present a collection of illustrative results that arise naturally from the analysis developed in the main text. While these are not used directly in the proofs or arguments of the preceding sections, they serve to highlight the scope and applicability of the main results and may offer insight for further exploration.
\begin{observation}
Let $n\in \mathbb N$ and $a,b,c,d\in \mathbb R^n\setminus \{0\}.$ Then $\{a\vee b, c\vee d\}$ is linearly dependent if and only if  at least
 one of the following statements is true.
\begin{enumerate}[wide=0pt]
 \item $\{a,c\}$ and $\{b,d\}$ are linearly dependent sets.
\item  $\{a,d\}$ and $\{b,c\}$ are linearly dependent sets.
\end{enumerate}
\end{observation}
\noindent Motivated by a result of \cite{bandyopadhyay2007differential} for the wedge product, we will see in the following example that for a subspace $ W\subseteq \Bvee^2(\mathbb R^2) $ with the property that $\dim W=2,$ $W\cap (x\vee \mathbb R^2)\neq \{0\}$ for all $x\in \mathbb R^2\setminus \{0\}$, but $W\neq \mathbb R^2\vee b$ for any $b\in \mathbb R^2\setminus \{0\}.$
\begin{example}\label{04/04/2024 example 1}
Suppose that $W=\spann \{e^1\vee (e^1+e^2), e^2\vee e^2\}.$ In such case, it is evident that $\dim W=2$ and that $W\cap (x\vee \mathbb R^2)\neq \{0\}$ for all $x\in \mathbb R^2\setminus \{0\}$, but $W\neq \mathbb R^2\vee b$ for any $b\in \mathbb R^2\setminus \{0\}.$ 
	
	
%
	
	
\end{example}
\begin{observation}
 For an $n$-dimensional subspace $W\subseteq \Bvee^2(\mathbb R^n)$, let us examine the next two claims: 
\begin{enumerate}[wide=0 pt]
\item  There exists $b\in \mathbb R^n \setminus \{0\}$ such that $W^{\perp}b=0.$
\item  $W\cap (x\vee \mathbb R^n)\neq \{0\}$ for all $x\in \mathbb R^n \setminus \{0\}.$
\end{enumerate}
 \Cref{17/11/2023 lemma 1} makes it clear that the assertion (1) will imply the statement (2). However, for $n=2$, we have one counterexample for (2) $ \nRightarrow$ (1) from \Cref{04/04/2024 example 1} and \Cref{17/11/2023 lemma 1} taken together.
\end{observation}
\begin{observation}
	In the situation of skew-symmetry, if we consider $W\subseteq \Lambda^2(\mathbb R^n)$ with $\dim W=n-1$, then the following three  statements are equivalent.
	\begin{enumerate}[wide= 0 pt]
		\item  There exists $b\in \mathbb R^n\setminus \{0\}$ such that $W^{\perp}b=0.$
		\item $W\cap (x\wedge \mathbb R^n)\neq \{0\}$ for all $x\in \mathbb R^n \setminus \{0\}.$
		\item   There exists $b\in \mathbb R^n\setminus \{0\}$ such that $W=\mathbb R^n \wedge b.$
	\end{enumerate}
\end{observation}
\begin{proof}
	Equivalence of the assertions  $(2)$ and $(3)$ follows from Theorem $4.6$ of \cite{bandyopadhyay2007differential}. To see $(1)$ implies $(3)$, let $A\in W^{\perp}$. Then $\langle A; x\wedge b\rangle=\langle b \lrcorner A; x\rangle=0$ for any $x\in \mathbb R^n.$ Hence $A\in (\mathbb R^n \wedge b)^{\perp}$. Therefore $W^{\perp}\subseteq (\mathbb R^n \wedge b)^{\perp}$ i.e., $\mathbb R^n \wedge b \subseteq W$ and as $\dim W=n-1,$ it follows that $W=\mathbb R^n \wedge b.$ To prove $(3)$ implies $(1)$, let $A\in W^{\perp}$. Then $\langle b\lrcorner A; x \rangle=\langle A; x\wedge b \rangle=0$ for all $x\in \mathbb R^n.$ Thus $b\lrcorner A=0$ for all $A\in W^{\perp}$. Hence $W^{\perp}b=0.$
\end{proof}

\noindent \textbf{Acknowledgement:} 
I am grateful to Prof. Saugata Bandyopadhyay for his guidance  in solving some problems  and for assisting me in substantiating many of the findings and writing this work.  
\section*{}

\end{document}